\documentclass[a4paper,11pt,reqno]{amsart}
\usepackage{amsmath}
\usepackage{amstext}
\usepackage{amsbsy}
\usepackage{amsopn}
\usepackage{upref}
\usepackage{amsthm}
\usepackage{amsfonts}
\usepackage{amssymb}
\usepackage{mathrsfs}
\usepackage{times}
\allowdisplaybreaks
 \newtheorem{theorem}{Theorem}

 \newtheorem{lemma}[theorem]{Lemma}
\theoremstyle{definition}
 
\theoremstyle{remark}
 
\newcommand{\ep}{\varepsilon}
\newcommand{\p}{\partial}

\begin{document}
\title[Schr\"odinger map]
      {Schr\"odinger flow into almost Hermitian manifolds}
\author[H.~Chihara]{Hiroyuki CHIHARA}
\address{Mathematical Institute,  
         Tohoku University, 
         Sendai 980-8578, Japan}
\email{chihara@math.tohoku.ac.jp}
\subjclass[2000]{Primary 53C44; Secondary 58J40, 47G30}
\thanks{Supported by JSPS Grant-in-Aid for Scientific Research \#20540151.}
\keywords{Schr\"odinger map, geometric analysis, pseudodifferential operators}
\begin{abstract}
We present a short-time existence theorem of solutions to 
the initial value problem for Schr\"odinger maps 
of a closed Riemannian manifold 
to a compact almost Hermitian manifold. 
The classical energy method cannot work for this problem 
since the almost complex structure of the target manifold is not supposed 
to be parallel with respect to the Levi-Civita connection. 
In other words, a loss of one derivative arises 
from the covariant derivative of the almost complex structure. 
To overcome this difficulty, 
we introduce a bounded pseudodifferential operator 
acting on sections of the pullback bundle, 
and essentially eliminate the loss of one derivative from 
the partial differential equation of the Schr\"odinger map. 
\end{abstract}
\maketitle
\section{Introduction}
\label{section:introduction}
Let $(M,g)$ be an $m$-dimensional closed Riemannian manifold 
with a Riemannian metric $g$, 
and let $(N,J,h)$ be a $2n$-dimensional compact almost Hermitian manifold 
with an almost complex structure $J$ and a Hermitian metric $h$. 
Consider the initial value problem for Schr\"odinger maps 
$u:\mathbb{R}{\times}M{\to}N$ of the form 
\begin{alignat}{2}
  \frac{\p{u}}{\p{t}}
& =
  J_u\tau(u)
& 
  \quad\text{in}\quad
& \mathbb{R}{\times}M,
\label{equation:pde}
\\
  u(0,x)
& =
  u_0(x)
& 
  \quad\text{in}\quad
& M,
\label{equation:data}
\end{alignat}
where 
$t\in\mathbb{R}$ is the time variable, 
$x{\in}M$, 
$\p{u}/\p{t}=du(\p/\p{t})$, 
$du$ is the differential of the mapping $u$, 
$u_0$ is a given map of $M$ to $N$, 
$\tau(u)=\operatorname{trace}\nabla{du}$ 
is the tension field of the map $u(t):{M}\to{N}$, 
and 
$\nabla$ is the induced connection. 
Here we observe local expression of $\tau(u)$. 
Let $x^1,\dotsc,x^m$ be local coordinates of $M$, 
and let $z^1,\dotsc,z^{2n}$ be local coordinates of $N$. 
We denote  
$$
g
=
\sum_{i,j=1}^m
g_{ij}dx^i{\otimes}dx^j, 
\quad
\sum_{k=1}^m
g_{ik}g^{kj}
=
\delta_{ij}, 
\quad
G=\det(g_{ij}), 
\quad
\Delta_g
=
\sum_{i,j=1}^m
\frac{1}{\sqrt{G}}
\frac{\p}{\p{x^i}}
g^{ij}\sqrt{G}\frac{\p}{\p{x^j}}, 
$$
$$
h=\sum_{a,b=1}^{2n}h_{ab}dz^a{\otimes}dz^b, 
\quad
\sum_{c=1}^{2n}
h_{ac}h^{cb}
=
\delta_{ab},
\quad
\Gamma^a_{bc}
=
\frac{1}{2}
\sum_{d=1}^{2n}
h^{ad}
\left(
\frac{\p{h_{bd}}}{\p{z^c}}
+
\frac{\p{h_{cd}}}{\p{z^b}}
-
\frac{\p{h_{bc}}}{\p{z^d}}
\right), 
$$
where $\delta_{ij}$ is Kronecker's delta. 
If we set $u^a=z^a{\circ}u$, 
the local expression of $\tau(u)$ is given by 
$$
\tau(u)
=
\sum_{i,j=1}^m
g^{ij}
\nabla{du}\left(\frac{\p}{\p{x^i}},\frac{\p}{\p{x^j}}\right)
=
\sum_{a=1}^{2n}
\left\{
\Delta_gu^a
+
\sum_{i,j=1}^m
\sum_{b,c=1}^{2n}
g^{ij}(x)
\Gamma^{a}_{bc}(u)
\frac{\p{u^b}}{\p{x^i}}
\frac{\p{u^c}}{\p{x^j}}
\right\} 
\left(\frac{\p}{\p{z^a}}\right)_u. 
$$
Then, \eqref{equation:pde} is a $2n\times2n$ system of 
quasilinear Schr\"odinger equations. 
\par
The equation \eqref{equation:pde} geometrically generalizes
two-sphere valued partial differential equations 
modeling the motion of vertex filament, 
ferromagnetic spin chain system and etc. 
See, e.g., \cite{darios}, \cite{hasimoto}, \cite{SSB} 
and references therein. 
In the last decade, these physics models have been generalized and
studied from a point of view of geometric analysis in mathematics. 
In other words, the relationship between 
the structure of the partial differential equation \eqref{equation:pde} 
and geometric settings have been investigated in the recent ten years. 
There are apparently two directions in the geometric analysis of partial
differential equations like \eqref{equation:pde}. 
\par
One of them is a geometric reduction of equations to 
simpler ones with values in the real or complex Euclidean space. 
This direction originated from Hasimoto's work \cite{hasimoto}. 
In their pioneering work \cite{CSU}, 
Chang, Shatah and Uhlenbeck first rigorously studied the PDE structure of 
\eqref{equation:pde} when $(M,g)$ is the real line with the usual metric 
and $(N,J,h)$ is a compact Riemann surface. 
They constructed a good moving frame along the map and 
reduced \eqref{equation:pde} to a simple complex-valued equation 
when $u(t,x)$ has a fixed base point as  $x\rightarrow+\infty$. 
Similarly, Onodera reduced a one-dimensional third or fourth order 
dispersive flow to a complex-valued equation in \cite{onodera2}. 
In \cite{NSU1} and \cite{NSU2}, 
Nahmod, Stefanov and Uhlenbeck obtained 
a system of semilinear Schr\"odinger equations 
from the equation of the Schr\"odinger map 
of the Euclidean space to the two-sphere 
when the Schr\"odinger map never takes values 
in some open set of the two-sphere. 
Nahmod, Shatah, Vega and Zeng constructed a moving frame along 
the Schr\"odinger map of the Euclidean space to a K\"ahler manifold 
in \cite{NSVZ}. 
Generally speaking, these reductions require some restrictions on 
the range of the mappings, and one cannot make use of them 
to solve the initial value problem for the original equations 
without restrictions on the range of the initial data. 
\par
The other direction of geometric analytic approach 
to partial differential equations like \eqref{equation:pde} 
is to consider how to solve the initial value problem. 
In his pioneering work \cite{koiso}, 
Koiso first reformulated the equation of 
the motion of vortex filament geometrically, 
and proposed the equation \eqref{equation:pde} 
when $(M,g)$ is the one-dimensional torus and 
$(N,J,h)$ is a compact K\"ahler manifold. 
Moreover, Koiso established the standard short-time existence theorem, 
and proved that if $(N,J,h)$ is locally symmetric, 
that is, $\nabla^NR=0$, then the solution exists globally in time, 
where $\nabla^N$ and $R$ are the Levi-Civita connection 
and the Riemannian curvature tensor of $N$ respectively. 
See \cite{PWW} also. 
Similarly, Onodera studied local and global existence theorem 
of solutions to a third-order dispersive flow for closed curves 
into K\"ahler manifolds in \cite{onodera1}. 
Gustafon, Kang and Tsai studied time-global stability 
of the Schr\"odinger map of the two-dimensional Euclidean space 
to the two-sphere around equivariant harmonic maps in \cite{GKT}. 
In \cite{DW}, Ding and Wang gave a short-time existence theorem of 
\eqref{equation:pde}-\eqref{equation:data} 
when $(M,g)$ is a general closed Riemannian manifold and 
$(N,J,h)$ is a compact K\"ahler manifold. 
However, they actually gave the proof only for the case that 
$(M,g)$ is the Euclidean space or the flat torus. 
We do not know whether their method of proof can work 
for a general closed Riemannian manifold $(M,g)$. 
Generally speaking, 
Schr\"odinger evolution equations are very delicate on lower order terms,  
in contrast with the heat equations, 
which can be easily treated together with any lower order term 
by the classical G{\aa}rding inequality. 
\par
Both of these two directions of geometric analysis of equations like
\eqref{equation:pde} are deeply concerned with the relationship between 
the geometric settings of equations and 
the theory of linear dispersive partial differential equations. 
For the latter subject, see, e.g., 
\cite{chihara2}, 
\cite{doi}, 
\cite{ichinose1}, 
\cite{ichinose2},  
\cite[Lecture VII]{mizohata} 
and references therein. 
Being concerned with the compactness of the source manifold, 
we need to mention local smoothing effect of dispersive partial
differential equations. 
It is well-known that solutions to the initial value problem for some
kinds of dispersive equations gain extra smoothness in comparison with
the initial data. 
In his celebrated work \cite{doi}, 
Doi characterized the existence of microlocal smoothing effect of
Schr\"odinger evolution equations on complete Riemannian manifolds 
according to the global behavior of the geodesic flow on the unit
cotangent sphere bundle over the source manifolds. 
Roughly speaking, 
the local smoothing effect occurs if and only if 
all the geodesics go to ``infinity''. 
In particular, if the source manifold is compact, 
then no smoothing effect occurs. 
For this reason, it is essential to study the initial value problem 
\eqref{equation:pde}-\eqref{equation:data} when $(M,g)$ is compact. 
\par
We should also mention the influence of the K\"ahler condition 
$\nabla^NJ=0$ on the structure of the equation \eqref{equation:pde}. 
All the preceding works on \eqref{equation:pde} assume that 
$(N,J,h)$ is a K\"ahler manifold. 
If $\nabla^NJ=0$, then \eqref{equation:pde} behaves like 
symmetric hyperbolic systems, 
and the classical energy method works effectively. 
See \cite{koiso} for the detail. 
If $\nabla^NJ\ne0$, then \eqref{equation:pde} has a first 
order terms in some sense, and the classical energy method breaks down. 
\par
The purpose of the present paper is to show 
a short-time existence theorem for
\eqref{equation:pde}-\eqref{equation:data} 
without the K\"ahler condition. 
To state our results, we here introduce function spaces of mappings.  
Set $\nabla_i=\nabla_{\p/\p{x^i}}$ for short. 
For a positive integer $k$, 
$H^k(M;TN)$ is the set of all continuous mappings 
$u:M{\rightarrow}N$ satisfying 
$$
\lVert{u}\rVert_{H^k}^2
=
\sum_{l=1}^k
\int_M
\lvert\nabla^lu\rvert^2
d\mu_g
<\infty,
$$
$$
\lvert\nabla^lu\rvert^2
=
\sum_{\substack{i_1,\dotsc,i_l \\ j_1,\dotsc,j_l=1}}^m
g^{i_1j_1}{\dotsb}g^{i_lj_l}
h
\left(
\nabla^lu\left(\frac{\p}{\p{x^{i_1}}},\dotsc,\frac{\p}{\p{x^{i_l}}}\right),
\nabla^lu\left(\frac{\p}{\p{x^{j_1}}},\dotsc,\frac{\p}{\p{x^{j_l}}}\right)
\right), 
$$
where $d\mu_g=\sqrt{G}dx^1{\dotsb}dx^m$ 
is the Riemannian measure of $(M,g)$. 
See e.g., \cite{hebey} for the Sobolev space of mappings.  
The Nash embedding theorem shows that 
there exists an isometric embedding 
$w{\in}C^\infty(N;\mathbb{R}^d)$ with some integer $d>2n$. 
See \cite{GR}, \cite{gunther} and \cite{nash}. 
Let $I$ be an interval in $\mathbb{R}$. 
We denote by $C(I;H^k(M;TN))$ the set of all 
$H^k(M;TN)$-valued continuous functions on $I$, 
In other words, we define it 
by the pullback of the function space by the isometry $w$ as 
$C(I;H^k(M;TN))=C(I;{w^\ast}H^k(M;\mathbb{R}^d))$, 
where $H^k(M;\mathbb{R}^d)$ is the usual Sobolev space 
of $\mathbb{R}^d$-valued functions on $M$. 
\par
Here we state our main results. 
\begin{theorem}
\label{theorem:main}
Let $k$ be a positive integer satisfying $2k>m/2+5$, 
and let $k_0$ be the minimum of $k$. 
Then, for any $u_0{\in}H^{2k}(M;TN)$, there exists 
$T=T(\lVert{u}\rVert_{H^{2k_0}})>$ such that 
{\rm \eqref{equation:pde}-\eqref{equation:data}} 
possesses a unique solution 
$u{\in}C([-T,T];H^{2k}(M;TN))$.  
\end{theorem}
Our strategy of proof consists of fourth order parabolic regularization 
and the uniform energy estimates of approximating solutions. 
Let $\Gamma(u^{-1}TN)$ be the set of sections of $u^{-1}TN$. 
Our idea of avoiding the difficulty due to $\nabla_iJ_u$ comes from 
diagonalization technique for some $2n{\times}2n$ system of
Schr\"odinger evolution equations developed in our work \cite{chihara1}. 
If we see \eqref{equation:pde} as a $2n{\times}2n$ system, 
$\nabla^NJ$ corresponds to some off-diagonal  blocks 
of the coefficient matrices of the first order terms. 
We introduce a bounded pseudodifferential operators acting on 
$\Gamma(u^{-1}TN)$ to eliminate $\nabla_iJ_u$ 
by using a transformation of u with this pseudodifferential operator. 
We here remark that 
$(\nabla_{\p/\p{z^a}}J)J=-J(\nabla_{\p/\p{z^a}}J)$ 
since $J^2=C^2_1J{\otimes}J=-I$, 
where $I$ is the identity mapping on $TN$ and 
$C^2_1$ is a contraction of $(2,2)$-tensor.  
This fact is the key to construct the pseudodifferential operator. 
We evaluate $\tilde{\Delta}_g^{l-1}\tau(u)$, ($l=1,\dotsc,k$), 
since 
$$
\tilde{\Delta}_g
=
\frac{1}{\sqrt{G}}
\sum_{i,j=1}^m
\nabla_ig^{ij}\sqrt{G}\nabla_j
$$
commutes with $\tau(u)$ and 
is never an obstruction to the energy estimates. 
For this reason, we need to use even order Sobolev space $H^{2k}(M;TN)$. 
It is easy to check that 
$\tilde{\Delta}_g$ is invariant under the change of
variables of $M$ and $N$.  
Indeed, one can check that $\tilde{\Delta}_g$ is invariant 
under the change of variables of $M$ in the same way as $\Delta_g$, 
and $\nabla_iV$ is invariant under the change of variables of $N$ 
for any section $V{\in}\Gamma(u^{-1}TN)$. 
Hence $\tilde{\Delta}_g$ is globally well-defined on $\Gamma(u^{-1}TN)$. 
\par
The plan of the present paper is as follows. 
Section~\ref{section:approximation} is devoted to parabolic regularization. 
In Section~\ref{section:proof} we prove Theorem~\ref{theorem:main}. 
\section{Parabolic Regularization}
\label{section:approximation}
This section is devoted to the short-time existence of solutions to 
the forward initial value problem 
for a fourth order parabolic equations of the form
\begin{alignat}{2}
  \frac{\p{u}}{\p{t}}
& =
  -\ep\tilde{\Delta}_g\tau(u)+J_u\tau(u)
& 
  \quad\text{in}\quad
& (0,\infty){\times}M,
\label{equation:pde-ep}
\\
  u(0,x)
& =
  u_0(x)
& 
  \quad\text{in}\quad
& M,
\label{equation:data-ep}
\end{alignat}
where $\ep\in(0,1]$ is a parameter. 
Roughly speaking, 
\eqref{equation:pde-ep}-\eqref{equation:data-ep} 
can be solved in the same way as the equation of harmonic heat flow 
$\p{u}/\p{t}=\tau(u)$. 
See e.g., \cite[Chapters~3 and 4]{nishikawa} 
for the study of the harmonic heat flow. 
In this section we shall show the following. 
\begin{lemma}
\label{theorem:approximation}
Let $l$ be an integer satisfying 
$l\geqslant{l_0}=[m/2]+4$. 
Then, for any $u_0{\in}H^l(M;TN)$, 
there exists 
$T_\ep=T(\ep,\lVert{u_0}\rVert_{H^{l_0}})>0$ 
such that the forward initial value problem 
{\rm \eqref{equation:pde-ep}-\eqref{equation:data-ep}} 
possesses a unique solution 
$u_\ep{\in}C([0,T_\ep];H^l(M;TN))$. 
\end{lemma}
We consider $\mathbb{R}^d$-valued equation pushed by $dw$. 
We split the proof of Lemma~\ref{theorem:approximation} into two steps. 
Firstly, we construct a solution 
taking values in a tubular neighborhood of $w(N)$. 
Secondly, we check that the solution is $w(N)$-valued. 
\par
If $u$ is a solution to \eqref{equation:pde-ep}, 
then $v=w{\circ}u$ satisfies 
$$
\frac{\p{v}}{\p{t}}
=
dw\left(\frac{\p{u}}{\p{t}}\right)
=
dw\left(-\ep\tilde{\Delta}_g\tau(u)+J_u\tau(u)\right)
=
-\ep\Delta_g^2v+F(v), 
$$
where $F(v)$ is of the form 
$$
F(v)
=
\ep
F_1(x,v,\bar{\nabla}v,\bar{\nabla}^2v,\bar{\nabla}^3v)
+
F_2(x,v,\bar{\nabla}v,\bar{\nabla}^2v)
$$
satisfying 
$F_1(x,y,0,0,0)=0$, $F_2(x,y,0,0)=0$ 
for any $(x,y){\in}M{\times}w(N)$, 
and $\bar{\nabla}$ is the connection induced by
$v(t):M\rightarrow{w(N)}$. 
\par
To solve this equation for $v$, 
we need elementary facts on the fundamental solution. 
Let $H^s(M)=(1-\Delta_g)^{-s/2}L^2(M)$ 
be the usual Sobolev space on $M$ of order $s\in\mathbb{R}$. 
We denote by $\mathscr{L}(\mathscr{H}_1,\mathscr{H}_j)$ 
the set of all bounded linear operators of 
a Hilbert space $\mathscr{H}_1$ to a Hilbert space $\mathscr{H}_2$. 
Set 
$\mathscr{L}(\mathscr{H}_1)=\mathscr{L}(\mathscr{H}_1,\mathscr{H}_1)$ 
for short. 
The existence and the properties of the fundamental solution 
are the following. 
\begin{lemma}
\label{theorem:psdo}
There exists an operator $E(t)$ satisfying 
$$
E(t)
\in
C((0,\infty);\mathscr{L}(H^s(M)))
\cap
C^1((0,\infty);\mathscr{L}(H^{s+4}(M),H^s(M)))
$$
for any $s\in\mathbb{R}$, such that 
$$
\left(
\frac{\p}{\p{t}}
+
\ep\Delta_g^2
\right)
E(t)
=
0
\quad\text{in}\quad
C((0,\infty);\mathscr{L}(H^{s+4}(M),H^s(M))), 
$$
$$
\lim_{t\downarrow0}E(t)v=v
\quad\text{for any}\quad
v{\in}H^s(M)
$$
$$
\left\{E(t)t^{3/4}\right\}_{t>0}
\quad\text{is bounded in}\quad
\mathscr{L}(H^{s-3}(M),H^s(M)).
$$
\end{lemma}
Lemma~\ref{theorem:psdo} is proved 
by the symbolic calculus of pseudodifferential operators. 
See \cite[Chapter~7, Theorem~4.1]{kumano-go} and \cite{iwasaki}. 
\par
Let $\delta>0$ be a sufficiently small constant. 
We denote by $w_\delta(N)$ 
a tubular neighborhood of $w(N)$ in $\mathbb{R}^d$, 
that is, 
$$
w_\delta(N)
=
\left\{
v_1+v_2\in\mathbb{R}^d
\ \vert \ 
v_1{\in}w(N),\ 
v_2{\in}T_{v_1}w(N)^{\perp},\ 
\lvert{v_2}\rvert<\delta
\right\},
$$
where 
$\lvert{v_2}\rvert=\sqrt{(v_2^1)^2+\dotsb+(v_2^d)^2}$ 
for 
$v_2=(v_2^1,\dotsc,v_2^d)\in\mathbb{R}^d$. 
Let $\pi:w_\delta(N){\rightarrow}w(N)$ be the projection 
defined by 
$\pi(v_1+v_2)=v_1$ 
for 
$v_1{\in}w(N)$, 
$v_2{\in}T_{v_1}w(N)^{\perp}$, 
$\lvert{v_2}\rvert<\delta$. 
We will solve the initial value problem 
\begin{equation}
\frac{\p{v}}{\p{t}}
=
-\ep\Delta_g^2v+F(\pi(v)), 
\quad
v(0,x)=w(u_0)(x),
\label{equation:aux}
\end{equation}
which is equivalent to an integral equation 
$v=\Phi(v)$, where 
$$
\Phi(v)(t)
=
E(t)w(u_0)
+
\int_0^t
E(t-s)
F(\pi(v(s)))
ds. 
$$
We apply the contraction mapping theorem to the integral equation 
in the framework 
$$
X_T^l
=
\left\{
v{\in}C([0,T];H^l(M;Tw_\delta(N)))
\ \vert\ 
\lVert{v(t)}\rVert_{H^l}
\leqslant
2\lVert{w(u_0)}\rVert_{H^l}
\quad\text{for}\quad 
t\in[0,T]
\right\}
$$
with some small $T>0$ determined below. 
By using Lemma~\ref{theorem:psdo} and the Sobolev embeddings, we have 
\begin{align*}
  \sup_{t\in[0,T]}
  \lVert\Phi(v)(t)-E(t)w(u_0)\rVert_{H^l}
& \leqslant
  C_lT^{1/4},
\\
  \sup_{t\in[0,T]}
  \lVert\Phi(v)(t)-\Phi(v^\prime)(t)\rVert_{H^l}
& \leqslant
  C_lT^{1/4}
  \sup_{t\in[0,T]}
  \lVert{v(t)-v^\prime(t)}\rVert_{H^l}
\end{align*}
for any $v,v^\prime{\in}X_T^l$, 
where $C_l$ is a positive constant depending on $\lVert{w(u_0)}\rVert_{H^l}$. 
Thus, $\Phi$ is a contraction mapping of $X_T^l$ to itself 
provided that $T$ is sufficiently small. 
The contraction mapping theorem shows the existence of a unique solution
to the integral equation. 
\par
Next we check that the solution 
$v{\in}C([0,T];H^l(M;Tw_\delta(N)))$ 
to \eqref{equation:aux} is $w(N)$-valued. 
Set $\rho(v)=v-\pi(v)$ for short. 
We remark that 
\begin{align*}
  \frac{\p{v}}{\p{t}}
& =
  -
  \ep\Delta_g^2v
  +
  F(\pi(v))
\\
& =
  -
  \ep\Delta_g^2\rho(v)
  -
  \ep\Delta_g^2\pi(v)
  +
  F(\pi(v))
\\
& =
  -
  \ep\Delta_g^2\rho(v)
  +
  dw
  \left(
  -\ep\tilde{\Delta}_g\tau(w^{-1}\circ\pi(v))
  +
  J_{w^{-1}\circ\pi(v)}\tau(w^{-1}\circ\pi(v))
  \right). 
\end{align*}
Since 
$\p\pi(v)/\p{t}\perp\rho(v)$ 
and 
$\rho(v){\in}T_{\pi(v)}w(N)^{\perp}$, 
we deduce
\begin{align*}
  \frac{d}{dt}
  \int_M
  \lvert\rho\rvert^2
  d\mu_g
& =
  2
  \int_M
  \left\langle\frac{\p{\rho(v)}}{\p{t}},\rho(v)\right\rangle
  d\mu_g
\\
& =
  2
  \int_M
  \left\langle\frac{\p{v}}{\p{t}},\rho(v)\right\rangle
  d\mu_g
\\
& =
  2
  \int_M
  \left\langle-\ep\Delta_g^2\rho(v)+dw(\dotsb),\rho(v)\right\rangle
  d\mu_g
\\
& =
  -2\ep
  \int_M
  \left\langle\Delta_g^2\rho(v),\rho(v)\right\rangle
  d\mu_g
\\
& =
  -2\ep
  \int_M
  \left\langle\Delta_g\rho(v),\Delta_g\rho(v)\right\rangle
  d\mu_g
  \leqslant
  0,
\end{align*}
where $\langle\cdot,\cdot\rangle$ is 
the standard inner product in $\mathbb{R}^d$. 
We conclude that $\rho(v)=0$ and therefore $v$ is $w(N)$-valued 
since $\rho(v(0))=\rho(w(u_0))=0$ 
This completes the proof of Lemma~\ref{theorem:approximation}. 
\section{Uniform Energy Estimates}
\label{section:proof}
The present section proves Theorem~\ref{theorem:main} 
in three steps: existence, uniqueness and recovery of continuity. 
Let $u_\ep{\in}C([0,T_\ep];H^{2k}(M;TN))$ be the unique solution to 
\eqref{equation:pde-ep}-\eqref{equation:data-ep}. 
\begin{proof}[Existence]
We shall show that there exists $T>0$ which is independent of $\ep\in(0,1]$,  
such that $\{u_\ep\}_{\ep\in(0,1]}$ is bounded in 
$L^\infty(0,T;H^{2k}(M;TN))$, 
which is the set of all $H^{2k}$-valued essentially bounded functions on $(0,T)$. 
If this is true, then the standard compactness argument shows that 
there exist $u$ and a subsequence $\{u_\ep\}$ such that 
\begin{align*}
  u_\ep \longrightarrow u
& \quad\text{in}\quad
  C([0,T];H^{2k-1}(M;TN)),
\\
  u_\ep \longrightarrow u
& \quad\text{in}\quad
  L^\infty(0,T;H^{2k}(M;TN))
  \quad\text{weakly star},
\end{align*}
as $\ep\downarrow0$, 
and $u$ solves \eqref{equation:pde}-\eqref{equation:data} 
and is $H^{2k}$-valued weakly continuous in time. 
\par
For $V\in\Gamma(u^{-1}TN)$, set 
$$
\lVert{V}\rVert^2
=
\int_Mh(V,V)d\mu_g
$$
for short. 
In view of the Sobolev embeddings, 
$\lVert{u}\rVert_{H^{2k}}$ is equivalent to  
$$
\left\{
\sum_{l=1}^k\lVert{\tilde{\Delta}_g^lu}\rVert^2
\right\}^{1/2}, 
\quad
\tilde{\Delta}_g^lu
=
\tilde{\Delta}_g^{l-1}\tau(u), 
$$
as a norm. 
We shall evaluate this for $u_\ep$ instead of $\lVert{u_\ep}\rVert_{H^{2k}}$. 
\par
The properties of the torsion tensor and the Riemannian curvature tensor 
show that for any vector fields $X$ and $Y$ on $M$, 
and for any $V\in\Gamma(u^{-1}TN)$, 
\begin{align}
  \nabla_Xdu(Y)
& =
  \nabla_Ydu(X)
  +
  du([X,Y]),
\label{equation:commutator1}
\\
  \nabla_X\nabla_YV
& =
  \nabla_Y\nabla_XV
  +
  \nabla_{[X,Y]}V
  +
  R(du(X),du(Y))V. 
\label{equation:commutator2}
\end{align}
Set 
$\nabla_t=\nabla_{\p/\p{t}}$, 
$X_i=\p/\p{x^i}$ 
and 
$Y_i=\sum_{j=1}^mg^{ij}X_j$ 
for short. 
\par
In what follows we express $u_\ep$ by $u$ simply. 
Any confusion will not occur. 
We sometimes write $\tilde{\Delta}_gu=\tau(u)$ and $Xu=du(X)$. 
We evaluate $\mathcal{N}_k(u)$ defined by 
$$
\mathcal{N}_k(u)^2
=
\sum_{j=1}^{k-1}
\lVert\tilde{\Delta}_g^lu\rVert^2
+
\lVert\Lambda\tilde{\Delta}_g^lu\rVert^2,
$$
where $\Lambda=\Lambda_\ep$ is a pseudodifferential operator defined later. 
Set 
$$
T_\ep^\ast
=
\sup\{
T>0
\ \vert \ 
\mathcal{N}_k(u(t))\leqslant2\mathcal{N}_k(u_0)
\ \text{for}\ 
t\in[0,T]
\}.
$$
To obtain the uniformly estimates of $\mathcal{N}_k(u)$, 
we need to compute
$$
\tilde{\Delta}_g^l
\left(
\frac{\p{u}}{\p{t}}
+
\ep\tilde{\Delta}_g^2u
-
J_u
\tilde{\Delta}_gu
\right)
=
0,
\quad
l=1,\dotsc,k.
$$
\par
In view of \eqref{equation:commutator1} and \eqref{equation:commutator2}, 
we have 
\begin{align*}
  \tilde{\Delta}_g\frac{\p{u}}{\p{t}} 
& =
  \frac{1}{\sqrt{G}}
  \sum_{i=1}
  \nabla_{X_i}
  \left\{
  \sqrt{G}\nabla_{Y_i}
  du(\p/\p{t})
  \right\}
\\
& =
  \frac{1}{\sqrt{G}}
  \sum_{i=1}
  \nabla_{X_i}
  \left\{
  \nabla_t\sqrt{G}du(Y_i)
  +
  \sqrt{G}
  du([Y_i,\p/\p{t}])
  \right\}
\\
& =
  \frac{1}{\sqrt{G}}
  \sum_{i=1}
  \nabla_{X_i}
  \left\{
  \nabla_t\sqrt{G}du(Y_i)
  \right\}
\\
& =
  \nabla_t
  \left\{
  \frac{1}{\sqrt{G}}
  \sum_{i=1}
  \nabla_{X_i}\sqrt{G}du(Y_i)
  \right\}
\\
& \quad
  +
  \frac{1}{\sqrt{G}}
  \sum_{i=1}
  \nabla_{[X_i,\p/\p{t}]}\sqrt{G}du(Y_i)
\\
& \quad
  +
  \sum_{i=1}^m
  R(du(X_i),du(\p/\p{t}))du(Y_i)
\\
& =
  \nabla_t\tilde{\Delta}_gu
  +
  \sum_{i=1}^m
  R(du(X_i),du(\p/\p{t}))du(Y_i). 
\end{align*}
Repeating this computation and using \eqref{equation:pde-ep}, we obtain 
\begin{align}
  \tilde{\Delta}_g^l
  \frac{\p{u}}{\p{t}}
& =
  \nabla_t\tilde{\Delta}_g^lu
  +
  \sum_{p=0}^{l-1}
  \tilde{\Delta}_g^{l-1-p}
  \left\{
  R(du(X_i),du(\p/\p{t}))\nabla_{Y_i}\tilde{\Delta}_g^pu
  \right\}
\nonumber
\\
& =
  \nabla_t\tilde{\Delta}_g^lu
  +
  \sum_{p=0}^{l-1}
  \tilde{\Delta}_g^{l-1-p}
  \left\{
  R(du(X_i),-\ep\tilde{\Delta}_g^2u+J_u\tilde{\Delta}_gu)
  \nabla_{Y_i}\tilde{\Delta}_g^pu
  \right\}
\nonumber
\\
& =
  \nabla_t\tilde{\Delta}_g^lu
  -
  \ep{P_{1,l}}
  -
  Q_{1,l},
\label{equation:aibu}
\\
  P_{1,l}
& =
  \sum_{p=0}^{l-1}
  \tilde{\Delta}_g^{l-1-p}
  \left\{
  R(du(X_i),\tilde{\Delta}^2u)
  \nabla_{Y_i}\tilde{\Delta}_g^pu
  \right\},
\nonumber
\\
  Q_{1,l}
& =
  -
  \sum_{p=0}^{l-1}
  \tilde{\Delta}_g^{l-1-p}
  \left\{
  R(du(X_i),J_u\tilde{\Delta}_gu)
  \nabla_{Y_i}\tilde{\Delta}_g^pu
  \right\}.
\nonumber
\end{align}
The Sobolev embeddings show that for $l=1,\dotsc,k$, 
there exists a constant $C_k>1$ depending only on $k\geqslant{k_0}$ and 
$\mathcal{N}_k(u_0)$ such that for $t\in[0,T_\ep^\ast]$ 
$$
\lVert{P_{1,l}}\rVert
\leqslant
C_k
\sum_{p=1}^{l+1}
\lVert\tilde{\Delta}_g^pu\rVert,
\quad
\lVert{Q_{1,l}}\rVert
\leqslant
C_k
\sum_{p=1}^l
\lVert\tilde{\Delta}_g^pu\rVert.
$$
Different positive constants depending only on $k\geqslant{k_0}$ and 
$\mathcal{N}_k(u_0)$ are denoted by the same notation $C_k$ below. 
\par
A direct computation shows that 
\begin{align}
  \tilde{\Delta}_g^l\left(J_u\tau(u)\right)
& =
  \frac{1}{\sqrt{G}}
  \sum_{i=1}^m
  \nabla_{X_i}
  \left(
  J_u\sqrt{G}\nabla_{Y_i}
  \tilde{\Delta}_g^lu
  \right)
\nonumber
\\
& \quad
  +
  l
  \sum_{i=1}^m
  (\nabla_{Y_i}J_u)\nabla_{X_i}
  \tilde{\Delta}_g^lu
\nonumber
\\
& \quad
  +
  (l-1)
  \sum_{i=1}^m
  (\nabla_{X_i}J_u)\nabla_{Y_i}
  \tilde{\Delta}_g^lu
  +
  Q_{2,l}
\nonumber
\\
& =
  \frac{1}{\sqrt{G}}
  \sum_{i,j=1}^m
  \nabla_i
  \left(
  g^{i,j}\sqrt{G}J_u\nabla_j
  \tilde{\Delta}_g^lu
  \right)
\nonumber
\\
& \quad
  +
  (2l-1)
  \sum_{i,j=1}^m
  g^{ij}
  (\nabla_iJ_u)\nabla_j
  \tilde{\Delta}_g^lu
  +
  Q_{2,l},
\label{equation:saki}
\end{align}
where $Q_{2,l}$ is a linear combination of terms of the form 
$(\nabla^{p+2}J_u)\nabla^{2l-p}u$, 
$p=0,1,\dotsc,2l-2$, 
and has the same estimate as $Q_{1,l}$. 
\par
Combining \eqref{equation:aibu} and \eqref{equation:saki}, we have 
\begin{align}
  \left\{
  \nabla_t
  +
  \ep\tilde{\Delta}_g^2
  -
  \frac{1}{\sqrt{G}}
  \sum_{i,j=1}^m
  \nabla_i
  g^{i,j}\sqrt{G}J_u\nabla_j
  -
  (2l-1)
  \sum_{i,j=1}^m
  g^{ij}
  (\nabla_iJ_u)\nabla_j
  \right\}
  \tilde{\Delta}_g^lu
& 
\nonumber
\\
  =
  {\ep}P_{1,l}+Q_{1,l}+Q_{2,l}.
&
\label{equation:pde2}
\end{align}
Set 
$$
P_{2,l}
=
(2l-1)
\sum_{i,j=1}^m
g^{ij}
(\nabla_iJ_u)\nabla_j
\tilde{\Delta}_g^lu
+ 
{\ep}P_{1,l}+Q_{1,l}+Q_{2,l}.  
$$
for short. 
$P_{2,l}$ can be estimated in the same way as $P_{1,l}$. 
Using \eqref{equation:pde2}, we deduce 
\begin{align}
  \frac{d}{dt}
  \sum_{l=1}^{k-1}
  \lVert\tilde{\Delta}_g^lu\rVert^2 
& =
  \frac{d}{dt}
  \sum_{l=1}^{k-1}
  \int_M
  h(\tilde{\Delta}_g^lu,\tilde{\Delta}_g^lu)
  d\mu_g
\nonumber
\\
& =
  2
  \sum_{l=1}^{k-1}
  \int_M
  h(\nabla_t\tilde{\Delta}_g^lu,\tilde{\Delta}_g^lu)
  d\mu_g
\nonumber
\\
& =
  -2\ep
  \sum_{l=1}^{k-1}
  \int_M
  h(\tilde{\Delta}_g^2\tilde{\Delta}_g^lu,\tilde{\Delta}_g^lu)
  d\mu_g
\nonumber
\\
& \quad
  +
  2
  \sum_{l=1}^{k-1}
  \sum_{i,j=1}^m
  \int_M
  \frac{1}{\sqrt{G}}
  h\left(
   \nabla_ig^{ij}\sqrt{G}J_u\nabla_j
   \tilde{\Delta}_g^lu,
   \tilde{\Delta}_g^lu
   \right)
  d\mu_g
\nonumber
\\
& \quad
  +
  2
  \sum_{l=1}^{k-1}
  \int_M
  h(P_{2,l},\tilde{\Delta}_g^lu)
  d\mu_g
\nonumber
\\
& =
  -2\ep
  \sum_{l=1}^{k-1}
  \int_M
  h(\tilde{\Delta}_g^{l+1}u,\tilde{\Delta}_g^{l+1}u)
  d\mu_g
\nonumber
\\
& \quad
  -
  2
  \sum_{l=1}^{k-1}
  \sum_{i,j=1}^m
  \int_M
  g^{ij}
  h(J_u\nabla_j\tilde{\Delta}_g^lu,\nabla_i\tilde{\Delta}_g^lu)
  d\mu_g
\nonumber
\\
& \quad
  +
  2
  \sum_{l=1}^{k-1}
  \int_M
  h(P_{2,l},\tilde{\Delta}_g^lu)
  d\mu_g
\nonumber
\\
& =
  -
  2\ep
  \sum_{l=1}^{k-1}
  \lVert\tilde{\Delta}_g^{l+1}u\rVert^2
  +
  2
  \sum_{l=1}^{k-1}
  \int_M
  h(P_{2,l},\tilde{\Delta}_g^lu)
  d\mu_g
\nonumber
\\
& \leqslant
  2C_k
  \sum_{l=1}^k
  \lVert\tilde{\Delta}_g^lu\rVert^2.
\label{equation:shiota}    
\end{align}
\par
To complete the energy estimates, 
we need to eliminate the first order term in \eqref{equation:pde2}. 
For this purpose, we here give the definition of 
the pseudodifferential operator $\Lambda$. 
Let $\{M_\nu\}$ be the finite set of local coordinate neighborhood of $M$, 
and let $x_\nu^1,\dotsb,x_\nu^m$ be the local coordinates in $M_\nu$. 
Let $\{N_\alpha\}$ be the set of local coordinate neighborhood of $N$, 
and let $z_\alpha^1,\dotsc,z_\alpha^{2n}$ be the local coordinates of $N_\alpha$. 
We denote by $C^\infty_0$ the set of all smooth functions with a compact support. 
Pick up partitions of unity 
$\{\phi_\nu\}{\subset}C^\infty_0(M)$ 
and 
$\{\Phi_\alpha\}{\subset}C^\infty_0(N)$ 
subordinated to $\{M_\nu\}$ and $\{N_\alpha\}$ respectively. 
Take $\{\psi_\nu\}, \{\xi_\nu\} \subset C^\infty_0(N)$ 
and $\{\Psi_\alpha\}, \{\Xi_\alpha\} \subset C^\infty_0(N)$ so that 
\begin{alignat*}{3}
  \xi_\nu=1
& \quad\text{in}\quad
  \operatorname{supp}[\phi_\nu],
& \quad
  \psi_\nu=1
& \quad\text{in}\quad
  \operatorname{supp}[\xi_\nu],
& \quad
  \operatorname{supp}[\psi_\nu]
& \subset
  M_\nu,
\\
  \Xi_\alpha=1
& \quad\text{in}\quad
  \operatorname{supp}[\Phi_\alpha],
& \quad
  \Psi_\alpha=1
& \quad\text{in}\quad
  \operatorname{supp}[\Xi_\alpha],
& \quad
  \operatorname{supp}[\Psi_\alpha]
& \subset
  N_\alpha.   
\end{alignat*}
We define a local $(1,1)$-tensor by 
$$
B_{\nu,\alpha,j}
=
-(2k-1)
\sum_{i=1}^m
g^{ij}(\nabla_iJ_u)
\quad\text{in}\quad
M_\nu{\cap}u^{-1}(N_\alpha).
$$
Note that $J_uB_{\nu,\alpha,j}=-B_{\nu,\alpha,j}J_u$. 
Using the partitions of unity and the local tensors, 
we define a properly supported pseudodifferential operators of order $-1$ 
acting on $\Gamma(u^{-1}TN)$ by 
\begin{align*}
  \tilde{\Lambda}
& =
  \frac{1}{2}
  \sum_\nu\sum_\alpha
  \phi_\nu(x)\Phi_\alpha(u)
  \sum_{j=1}^m
  J_uB_{\nu,\alpha,j}\nabla_j
  \xi_\nu(x)\Xi_\alpha(u)
  (1-\Delta_g)^{-1}
  \psi_\nu(x)\Psi_\alpha(u)
\\
& -
  \frac{1}{2}
  \sum_\nu\sum_\alpha
  \phi_\nu(x)\Phi_\alpha(u)
  \sum_{j=1}^m
  B_{\nu,\alpha,j}J_u\nabla_j
  \xi_\nu(x)\Xi_\alpha(u)
  (1-\Delta_g)^{-1}
  \psi_\nu(x)\Psi_\alpha(u),  
\end{align*}
where $\nabla_j$ is expressed by the coordinates 
$x_\nu^1,\dotsc,x_\nu^m$ 
and 
$z_\alpha^1,\dotsc,z_\alpha^{2n}$. 
Here we remark that for 
$$
V
=
\sum_{a=1}^{2n}
V_{\nu,\alpha}^a
\left(\frac{\p}{\p{z_\alpha^a}}\right)_u
$$
supported in $M_\nu{\cap}u^{-1}(N_\alpha)$, 
$$
\xi_\nu(x)\Xi_\alpha(u)(1-\Delta_g)^{-1}V
=
\sum_{a=1}^{2n}
\left(
\xi_\nu(x)\Xi_\alpha(u)(1-\Delta_g)^{-1}V_{\nu,\alpha}^a
\right)
\left(\frac{\p}{\p{z_\alpha^a}}\right)_u
$$
is well-defined and supported in $M_\nu{\cap}u^{-1}(N_\alpha)$ also. 
We make use of elementary theory of pseudodifferential operators freely. 
We remark that 
each term in $\tilde{\Lambda}$ is properly supported 
and invariant under the change of coordinates in $M$ and $N$ 
up to cut-off functions. 
Then, we can deal with each term as if it were 
a pseudodifferential operator on $\mathbb{R}^m$ 
acting on $\mathbb{R}^d$-valued functions. 
Then, we can make use of pseudodifferential operators 
whose symbols have limited smoothness. 
See \cite[Section~2]{chihara1} and \cite{nagase}. 
In other words, we do not have to take care of the type of 
this pseudodifferential operator. 
It is well-known that the type of pseudodifferential operators 
on manifolds have some restrictions in general. 
\par
Set $\Lambda=I-\tilde{\Lambda}$ and $\Lambda^\prime=I+\tilde{\Lambda}$, 
where $I$ is the identity mapping. 
Since $\Lambda^\prime\Lambda=I-\tilde{\Lambda}^2$ 
and $\tilde{\Lambda}^2$ is a pseudodifferential operator of order $-2$, 
we deduce that for $t\in[0,T_\ep^\ast]$, 
$$
C_k^{-1}\mathcal{N}_k(u)^2
\leqslant
\sum_{l=1}^k
\lVert\tilde{\Delta}_g^lu\rVert^2
\leqslant
C_k\mathcal{N}_k(u)^2.
$$
\par
We apply $\Lambda$ to \eqref{equation:pde2} with $l=k$. 
A direct computation shows that 
$$
\Lambda\nabla_t=\nabla_t\Lambda-\frac{\p\Lambda}{\p{t}},
\quad
\left\lVert
\frac{\p\Lambda}{\p{t}}
\tilde{\Delta}_g^ku
\right\rVert
\leqslant
C_k\mathcal{N}_k(u)
$$
for $t\in[0,T_\ep^\ast]$. 
The matrices of principal symbols of $\Lambda$ and $\tilde{\Delta}_g^2$ 
commute with each other since the matrix of the principal symbol of 
$\tilde{\Delta}_g^2$ is $(g^{ij}\xi_i\xi_j)^2I_{2n}$, 
where $I_{2n}$ is the $2n\times2n$ identity matrix. 
Then, 
$[\Lambda,\tilde{\Delta}_g^2]=[\tilde{\Lambda},\tilde{\Delta}_g^2]$ 
is a pseudodifferential operator of order $2$. 
Hence we have 
$$
\Lambda\tilde{\Delta}_g^2
=
\tilde{\Delta}_g^2\Lambda
+
[\Lambda,\tilde{\Delta}_g^2], 
\quad
\lVert[\Lambda,\tilde{\Delta}_g^2]\tilde{\Delta}_g^ku\rVert
\leqslant
C_k\lVert\tilde{\Delta}_g^{k+1}u\rVert
$$
for $t\in[0,T_\ep^\ast]$. 
\par
Here we set 
$$
A
=
\sum_{i,j=1}^m
\frac{1}{\sqrt{G}}
\nabla_ig^{ij}\sqrt{G}J_u\nabla_j
$$
for short. 
Since $I=\Lambda^\prime\Lambda+\tilde{\Lambda}^2$, 
we deduce that 
$$
\Lambda(-A)
=
-\Lambda{A}(\Lambda^\prime\Lambda+\tilde{\Lambda}^2)
=
-A+\tilde{\Lambda}A-A\tilde{\Lambda}
+(\tilde{\Lambda}A\tilde{\Lambda}+A\tilde{\Lambda}^2),
$$
and $\tilde{\Lambda}A\tilde{\Lambda}+A\tilde{\Lambda}^2$ 
is $L^2$-bounded for $t\in[0,T_\ep^\ast]$. 
The principal symbol of $(1-\Delta_g)^{-1}$ is globally defined as 
$I_{2n}/g^{ij}\xi_i\xi_j$. 
We deduce that modulo $L^2$-bounded operators, 
\begin{align*}
  \tilde{\Lambda}A
& \equiv
  -
  \frac{1}{2}
  \sum_\nu\sum_\alpha
  \phi_\nu(x)\Phi_\alpha(u)
  \sum_{j=1}^m
  B_{\nu,\alpha,j}J_u^2\nabla_j
  \xi_\nu(x)\Xi_\alpha(u)
  (1-\Delta_g)^{-1}
  \psi_\nu(x)\Psi_\alpha(u)
  \tilde{\Delta}_g
\\
& \equiv
  -
  \frac{1}{2}
  \sum_\nu\sum_\alpha
  \phi_\nu(x)\Phi_\alpha(u)
  \sum_{j=1}^m
  B_{\nu,\alpha,j}\nabla_j
  \xi_\nu(x)\Xi_\alpha(u)
  \psi_\nu(x)\Psi_\alpha(u)
\\
& =
  -
  \frac{1}{2}
  \sum_\nu\sum_\alpha
  \phi_\nu(x)\Phi_\alpha(u)
  \sum_{j=1}^m
  B_{\nu,\alpha,j}\nabla_j
\\
& =
  \frac{(2k-1)}{2}
  \sum_{i,j}g^{ij}(\nabla_iJ_u)\nabla_j,
\\
  -A\tilde{\Lambda}
& \equiv
  -
  \frac{1}{2}
  \sum_\nu\sum_\alpha
  \tilde{\Delta}_g
  \phi_\nu(x)\Phi_\alpha(u)
  \sum_{j=1}^m
  J_u^2B_{\nu,\alpha,j}\nabla_j
  \xi_\nu(x)\Xi_\alpha(u)
  (1-\Delta_g)^{-1}
  \psi_\nu(x)\Psi_\alpha(u)
\\
& \equiv
  -
  \frac{1}{2}
  \sum_\nu\sum_\alpha
  \phi_\nu(x)\Phi_\alpha(u)
  \sum_{j=1}^m
  B_{\nu,\alpha,j}\nabla_j
  \xi_\nu(x)\Xi_\alpha(u)
  \psi_\nu(x)\Psi_\alpha(u)
\\
& =
  -
  \frac{1}{2}
  \sum_\nu\sum_\alpha
  \phi_\nu(x)\Phi_\alpha(u)
  \sum_{j=1}^m
  B_{\nu,\alpha,j}\nabla_j
\\
& =
  \frac{(2k-1)}{2}
  \sum_{i,j}g^{ij}(\nabla_iJ_u)\nabla_j. 
\end{align*}
\par
Combining computations above, we obtain 
$$
\left\{
\nabla_t
+
\ep\tilde{\Delta}_g^2
-
\frac{1}{\sqrt{G}}
\sum_{i,j=1}^m
\nabla_i
g^{i,j}\sqrt{G}J_u\nabla_j
\right\}
\Lambda
\tilde{\Delta}_g^ku
=
{\ep}P_k+Q_k,
$$
where $P_k$ and $Q_k$ are estimated as 
$$
\lVert{P_k}\rVert
\leqslant
C_k(\lVert\tilde{\Delta}_g\Lambda\tilde{\Delta}_g^ku\rVert+\mathcal{N}_k(u)), 
\quad
\lVert{Q_k}\rVert
\leqslant
C_k\mathcal{N}_k(u).
$$
In the same way as \eqref{equation:shiota}, we deduce 
\begin{align}
  \frac{d}{dt}\lVert\Lambda\tilde{\Delta}_g^ku\rVert^2
& \leqslant
  -2\ep
  \lVert\tilde{\Delta}_g\Lambda\tilde{\Delta}_g^ku\rVert^2
  +
  C_k\ep
  \lVert\tilde{\Delta}_g\Lambda\tilde{\Delta}_g^ku\rVert
  \lVert\Lambda\tilde{\Delta}_g^ku\rVert
  +
  C_k
  \mathcal{N}_k(u)
  \lVert\Lambda\tilde{\Delta}_g^ku\rVert
\\
& \leqslant
  C_k
  \mathcal{N}_k(u)
  \lVert\Lambda\tilde{\Delta}_g^ku\rVert.
\label{equation:reiko}
\end{align}
Combining \eqref{equation:shiota} and \eqref{equation:reiko}, we obtain
\begin{equation}
\frac{d}{dt}\mathcal{N}_k(u)
\leqslant
C_k
\mathcal{N}_k(u) 
\quad\text{for}\quad
t\in[0,T_\ep^\ast].
\label{equation:joe}
\end{equation}
If we take $t=T_\ep^\ast$, then we have 
$2\mathcal{N}_k(u_0)\leqslant\mathcal{N}_k(u_0)e^{C_kT_\ep^\ast}$, 
which implies that 
$T_\ep^\ast{\geqslant}T=\log2/C_k>0$. 
Thus $\{u_\ep\}_{\ep{\in(0,1]}}$ is bounded in $L^\infty(0,T;H^{2k}(M;TN))$. 
This completes the proof. 
\end{proof}
\begin{proof}[Uniqueness]
Let $u_1,u_2{\in}L^\infty(0,T;H^{2k}(M;TN))$ be solutions to 
\eqref{equation:pde}-\eqref{equation:data}. 
Set $v_1=w{\circ}u_1$ and $v_2=w{\circ}u_2$ for short. 
We denote by $\Pi_v:\mathbb{R}^d \rightarrow \mathbb{R}^d$ 
the projection to $T_vw(N)$ for $v{\in}w(N)$. 
Since 
$$
\frac{\p{v}}{\p{t}}
=
\tilde{J}_{v_j}\Pi_{v_j}\Delta_gv_j,
\quad
\tilde{J}_{v_j}
=
du_j{\circ}J_{u_j}{\circ}du_j^{-1},
\quad
j=1,2, 
$$
$v=v_1-v_2$ solves 
\begin{align*}
  \frac{\p{v}}{\p{t}}
& =
  \tilde{J}_{v_1}\Pi_{v_1}\Delta_gv
  +
  \left(\tilde{J}_{v_1}\Pi_{v_1}-\tilde{J}_{v_2}\Pi_{v_2}\right)\Delta_gv_2
\\
& =
  \tilde{J}_{v_1}\Pi_{v_1}\Delta_gv
  +
  A(v_1,v_2,\Delta_gv_2)v. 
\end{align*}
Here we applied the mean value theorem to the second term of 
the right hand side of the above equation, 
and $A(v_1,v_2,\Delta_gv_2)$ is an appropriate $d{\times}d$ matrix. 
\par
Using the integration by parts, we have 
\begin{align}
  \frac{d}{dt}
  \int_M\langle{v,v}\rangle
  d\mu_g
& =
  2
  \int_M
  \left\langle\frac{\p{v}}{\p{t}},v\right\rangle
  d\mu_g
\nonumber
\\
& =
  2
  \int_M
  \left\langle{
  \tilde{J}_{v_1}\Pi_{v_1}\Delta_gv
  +
  A(v_1,v_2,\Delta_gv_2)v,
  v
  }\right\rangle
  d\mu_g
\nonumber
\\
& \leqslant
  C
  \int_M
  \left\{
  \langle{v,v}\rangle
  +
  \sum_{i,j=1}^m
  g^{ij}
  \left\langle
  \frac{\p{v}}{\p{x^i}},
  \frac{\p{v}}{\p{x^j}}
  \right\rangle
  \right\}
  d\mu_g.
\label{equation:katase}
\end{align}
Using the properties of $\tilde{J}_{v_1}$ and the projection, 
and the integration by parts, we deduce 
\begin{align}
  \frac{d}{dt}
  \int_M
  \sum_{i,j=1}^m
  g^{ij}
  \left\langle
  \frac{\p{v}}{\p{x^i}},
  \frac{\p{v}}{\p{x^j}}
  \right\rangle
  d\mu_g
& =
  2
  \int_M
  \sum_{i,j=1}^m
  g^{ij}
  \left\langle
  \frac{\p^2{v}}{\p{t}\p{x^i}},
  \frac{\p{v}}{\p{x^j}}
  \right\rangle
  d\mu_g
\nonumber
\\
& =
  -2
  \int_M
  \left\langle
  \frac{\p{v}}{\p{t}},
  \Delta_gv
  \right\rangle
  d\mu_g
\nonumber
\\
& =
  -
  2
  \int_M
  \left\langle
  \tilde{J}_{v_1}\Pi_{v_1}\Delta_gv+Av,
  \Delta_gv
  \right\rangle
  d\mu_g
\nonumber
\\
& =
  -
  2
  \int_M
  \left\langle
  \tilde{J}_{v_1}\Pi_{v_1}\Delta_gv,
  \Pi_{v_1}\Delta_gv
  \right\rangle
  d\mu_g
\nonumber
\\
& \quad 
  -
  2
  \int_M
  \langle{Av,\Delta_gv}\rangle
  d\mu_g
\nonumber
\\
& =
  -
  2
  \int_M
  \langle{Av,\Delta_gv}\rangle
  d\mu_g
\nonumber
\\
& \leqslant
  C
  \int_M
  \left\{
  \langle{v,v}\rangle
  +
  \sum_{i,j=1}^m
  g^{ij}
  \left\langle
  \frac{\p{v}}{\p{x^i}},
  \frac{\p{v}}{\p{x^j}}
  \right\rangle
  \right\}
  d\mu_g.
\label{equation:nana}
\end{align}
Combining \eqref{equation:katase} and \eqref{equation:nana}, we have 
$$
\frac{d}{dt}
\int_M
\left\{
\langle{v,v}\rangle
+
\sum_{i,j=1}^m
g^{ij}
\left\langle
\frac{\p{v}}{\p{x^i}},
\frac{\p{v}}{\p{x^j}}
\right\rangle
\right\}
d\mu_g
\leqslant
C
\int_M
\left\{
\langle{v,v}\rangle
+
\sum_{i,j=1}^m
g^{ij}
\left\langle
\frac{\p{v}}{\p{x^i}},
\frac{\p{v}}{\p{x^j}}
\right\rangle
\right\}
d\mu_g,
$$
which implies $v=0$. 
This completes the proof. 
\end{proof}
\begin{proof}[Continuity in time]
Let $u{\in}L^\infty(0,T;H^{2k}(M;TN))$ 
be the unique solution to \eqref{equation:pde}-\eqref{equation:data}. 
We remark that 
$u{\in}C([0,T];H^{2k-1}(M;TN))$ 
and $\tilde{\Delta}_g^ku$ is a weakly continuous 
$L^2(M;TN)$-valued function on $[0,T]$. 
We identify $N$ and $w(N)$ below. 
Let $\{u_\ep\}_{\ep\in(0,1]}$ be a sequence of solutions 
to \eqref{equation:pde-ep}-\eqref{equation:data-ep}, 
which approximates $u$. 
We can easily check that for any $\phi{\in}C^\infty([0,T]\times{M};\mathbb{R}^d)$, 
\begin{align*}
  \Lambda_\ep\phi \longrightarrow \Lambda\phi
& \quad\text{in}\quad
  L^2((0,T){\times}M;\mathbb{R}^d),
\\
  u_\ep \longrightarrow u
& \quad\text{in}\quad
  L^2((0,T){\times}M;\mathbb{R}^d),
\\
  \Lambda_\ep\tilde{\Delta}_g^ku_\ep \longrightarrow \tilde{u}
& \quad\text{in}\quad
  L^2((0,T){\times}M;\mathbb{R}^d)
  \quad\text{weakly star},
\end{align*}
as $\ep\downarrow0$ 
with some $\tilde{u}$. 
Then, $\tilde{u}=\Lambda\tilde{\Delta}_g^ku$ 
in the sense of distributions. 
The time-continuity of $\tilde{\Delta}_g^ku$ is equivalent to 
that of $\Lambda\tilde{\Delta}_g^ku$ since 
$\Lambda{\in}C([0,T];\mathscr{L}(L^2(M;\mathbb{R}^d)))$, 

\par
It suffices to show that  
\begin{equation}
\lim_{t\downarrow0}
\Lambda(t)\tilde{\Delta}_g^ku(t)
=
\Lambda(0)\tilde{\Delta}_g^ku_0
\quad\text{in}\quad
L^2(M;\mathbb{R}^d),
\label{equation:gillian}  
\end{equation}
since the other cases can be proved in the same way. 
\eqref{equation:joe} and the lower semicontinuity of $L^2$-norm imply 
$$
\sum_{l=1}^{k-1}
\lVert{\tilde{\Delta}_g^lu(t)}\rVert^2
+
\lVert{\Lambda(t)\tilde{\Delta}_g^ku(t)}\rVert^2
\leqslant
\sum_{l=1}^{k-1}
\lVert{\tilde{\Delta}_g^lu_0}\rVert^2
+
\lVert{\Lambda(0)\tilde{\Delta}_g^ku_0}\rVert^2
+
C_k\mathcal{N}_k(u_0)^2t.
$$
Letting $t\downarrow0$, we have 
$$
\limsup_{t\downarrow0}
\lVert{\Lambda(t)\tilde{\Delta}_g^ku(t)}\rVert^2
\leqslant
\lVert{\Lambda(0)\tilde{\Delta}_g^ku_0}\rVert^2,
$$
which implies \eqref{equation:gillian}. 
This completes the proof.
\end{proof}
%
%

\end{document}